\documentclass{article}
\begin{document}
\newtheorem{proposition}{Proposition}[section]
\newtheorem{definition}{Definition}[section]
\newtheorem{lemma}{Lemma}[section]
 
\title{\bf Generalizations of Jordan Algebras and Malcev Algebras}
\author{Keqin Liu\\Department of Mathematics\\The University of British Columbia\\Vancouver, BC\\
Canada, V6T 1Z2}
\date{April, 2006}
\maketitle

\begin{abstract}We introduce two classes of nonassociative algebras and define the building blocks in the context of the new nonassociative algebras.
\end{abstract}

My generalizations of Jordan Algebras and Malcev Algebras come from $\mathcal{Z}_2$-algebras which lead to my attention in exploring the curiosity to generalize Lie algebras. The definitions of the generalized Jordan Algebras and the generalized Malcev Algebras are based on two passages. One passage is from associative $\mathcal{Z}_2$-algebras to generalized Jordan Algebras. The other passage is from alternative $\mathcal{Z}_2$-algebras to generalized Malcev Algebras. 

\medskip
In this paper, all vector spaces are vector spaces over fields of characteristic not two and three, and all associative algebras have an identity. 

\section{Generalized Jordan Algebras}

We begin this section with the following definition.

\begin{definition}\label{def1} Let $A=A_0\oplus A_1$ (as vector spaces) be an algebra.
\begin{description}
\item[(i)] $A$ is called a {\bf $\mathcal{Z}_2$-algebra} if 
$$A_0A_0\subseteq A_0,\quad A_0A_1+A_1A_0\subseteq A_1 \quad\mbox{and}\quad A_1A_1=0.$$
\item[(ii)] $A$ is called an {\bf associative $\mathcal{Z}_2$-algebra} if $A$ is both a $\mathcal{Z}_2$-algebra and an associative algebra.
\item[(iii)] $A$ is called an {\bf alternative $\mathcal{Z}_2$-algebra} if $A$ is both a $\mathcal{Z}_2$-algebra and an alternative algebra.
\end{description}
\end{definition}

The algebra $\mathcal{O}^{(2)}$ introduced in \cite{hl2} is a $16$-dimensional alternative $\mathcal{Z}_2$-algebra which is not associative.

\medskip
Using the one-sided commutative law introduced in \cite{Liu}, we now define generalized Jordan algebras.

\begin{definition}\label{def 2} A vector space $J$ is called a {\bf generalized Jordan algebra} if exists a binary operation $\bullet$: $J\times J \to J$ such that the following two properties hold.
\begin{description}
\item[(i)] The binary operation $\bullet$ is {\bf right commutative}; that is,
$$x\bullet (y\bullet z)=x\bullet (z\bullet y)\quad\mbox{for $x$, $y$, $z\in J$}.$$
\item[(ii)] The binary operation $\bullet$ satisfies the {\bf Jordan identity}:
$$
(y\bullet x)\bullet  (x\bullet x)=\left(y\bullet (x\bullet x)\right)\bullet x\quad\mbox{for $x$, $y\in J$}$$
and the {\bf Hu-Liu bullet identity}:
$$
x\bullet \left(y\bullet  (x\bullet x)\right)-(x\bullet y)\bullet (x\bullet x)
=2(x\bullet x)\bullet (y\bullet x)-2\left((x\bullet x)\bullet y\right)\bullet x
$$
for $x$, $y\in J$.
\end{description}
\end{definition}

A generalized Jordan algebra $J$ is also denoted by $(J,\, +,\, \bullet )$, where the binary operation $\bullet$ is called the {\bf bullet product}. A generalized Jordan algebra 
$(J,\, +,\, \bullet )$ is said to be {\bf right unital} if there exists an element $1$ of such that
\begin{equation}\label{eq1}
x\bullet 1=x \quad\mbox{for $x\in J$}.
\end{equation}
An element $1$ satisfying (\ref{eq1}) is called a {\bf right unit}.

\medskip
The following proposition establishes the passage from a associative $\mathcal{Z}_2$-algebra to a 
generalized Jordan algebra.

\begin{proposition}\label{pr1} If $A=A_0\oplus A_1$ is an associative $\mathcal{Z}_2$-algebra, then $(A,\, +,\, \bullet )$ is a right unital generalized Jordan algebra, where the bullet product $\bullet$ is defined by
\begin{equation}\label{eq2}
x\bullet y: =\frac12 (xy_0+y_0x)
\end{equation}
for $x\in A$, $y=y_0+y_1\in A$ and $y_i\in A_i$ with $i=0$ and $1$.
\end{proposition}

\medskip
\noindent
{\bf Proof} It is a straightforward computation.

\hfill\raisebox{1mm}{\framebox[2mm]{}}

\bigskip
\begin{definition}\label{def3} Let $I$ be a subspace of a generalized Jordan algebra 
$(J,\, +,\, \bullet )$.
\begin{description}
\item[(i)] $I$ is called an {\bf ideal} of $J$ if $I\bullet J\subseteq I$ and 
$J\bullet I\subseteq I$.
\item[(ii)] $I$ is called a {\bf (generalized Jordan) subalgebra} of $J$ if 
$I\bullet I\subseteq I$.
\end{description}
\end{definition}

\bigskip
The {\bf annihilator} $J^{ann}$ of a generalized Jordan algebra $(J,\, +,\, \bullet )$ over a field $\mathbf{k}$ is defined by
\begin{equation}\label{eq3}
J^{ann}:=\sum _{x, y\in J}\mathbf{k} (x\bullet y - y\bullet x).
\end{equation}
It is easy to check that $J^{ann}$ is an ideal of $J$ and $\displaystyle\frac{J}{J^{ann}}$ is a Jordan algebra.

\medskip
The building blocks in the context of generalized Jordan  algebras are defined by using annihilators.

\bigskip
\begin{definition}\label{def4} A generalized Jordan algebra $J$ is said to be {\bf simple} if $J^{ann}\ne 0$ and $J$ has no ideals which are not equal to $\{0\}$, $J^{ann}$ and $J$.
\end{definition}

If $(J,\, +,\, \bullet )$ is a right unital generalized Jordan algebra, then
$$ J^{ann}=\{\, a \, |\, 1\bullet a=0 \,\}$$
and
$$ \{\, 1+a \,|\, a\in J^{ann} \,\}=\mbox{the set of all right units of $J$}, $$
where $1$ is a right unit of $J$.

\medskip
The next proposition gives a basic property of annihilators.

\begin{proposition}\label{pr2} If $(J,\, +,\, \bullet )$ is a generalized Jordan algebra over a field $\mathbf{k}$, then the annihilator $J^{ann}$ becomes a bimodule over the Jordan algebra $\displaystyle\frac{J}{J^{ann}}$ under the bimodule action:
$$
\bar{x}\, u=u\, \bar{x}:=u\bullet x \quad\mbox{for $\bar{x}=x+ J^{ann}\in \displaystyle\frac{J}{J^{ann}}$,
$x\in J$ and $u\in J^{ann}$.}
$$
\end{proposition}

\medskip
\noindent
{\bf Proof} This is a consequence of the Hu-Liu bullet identity. See page 82 of \cite{J3} for the definition of bimodules over Jordan algebras.

\hfill\raisebox{1mm}{\framebox[2mm]{}}

\bigskip
Let $J$ and $\bar{J}$ be two generalized Jordan algebras. A linear map $\phi$: $J\to \bar{J}$ is called a {\bf homomorphism} if 
$$\phi (x\bullet y )=\phi(x)\bullet \phi(y) \quad\mbox{for $x$, $y\in J$.}$$

\begin{definition}\label{def5} A generalized Jordan algebra $J$ is said to be {\bf $\mathcal{Z}_2$-special} if there exists an associative $\mathcal{Z}_2$-algebra $A=A_0\oplus A_1$ such that there is an injective homomorphism from $J$ the generalized Jordan algebra $(A,\, +,\, \bullet)$, where the bullet product $\bullet$ in $A$ is defined by (\ref{eq2}). Generalized Jordan algebras which are not $\mathcal{Z}_2$-special are said be 
{\bf $\mathcal{Z}_2$-exceptional}.
\end{definition}

\section{Generalized Malcev Algebras}

We begin this section with the definition of generalized Malcev algebras.

\begin{definition}\label{def 2.1} A vector space $M$ is called a {\bf generalized Malcev algebra} if exists a binary operation $\langle \,, \,\rangle$: $M\times M \to M$ such that the following two properties hold.
\begin{description}
\item[(i)] The binary operation $\langle \,, \,\rangle$ is {\bf right anti-commutative}; that is,
$$\langle x, \langle y, z\rangle\rangle=-\langle x, \langle z, y\rangle\rangle\quad\mbox{for $x$, $y$, $z\in J$}.$$
\item[(ii)] The binary operation $\langle \,, \,\rangle$ satisfies the {\bf Hu-Liu angle identities}:
\begin{eqnarray}
&&\langle \langle x, y\rangle , \langle x, z\rangle\rangle -
\langle \langle x, z\rangle , \langle x, y\rangle\rangle \nonumber\\
\label{eq4} &=&\langle\langle\langle x, y\rangle , z\rangle , x\rangle -
\langle x, \langle\langle x, y\rangle , z\rangle \rangle +
\langle x, \langle x, \langle y, z\rangle \rangle \rangle -
\langle\langle x, \langle y, z\rangle \rangle , x\rangle +\nonumber\\
&&\qquad +\langle\langle x, \langle x, z\rangle \rangle , y\rangle -
\langle\langle \langle x, z\rangle , x\rangle , y\rangle ,
\end{eqnarray}
\begin{equation}\label{eq5}
\langle \langle y, x\rangle , \langle x, z\rangle\rangle =
\langle\langle\langle y, x\rangle , z\rangle , x\rangle -
\langle y, \langle\langle x, z\rangle , x\rangle \rangle - 
\langle\langle\langle y, z\rangle , x\rangle , x\rangle ,
\end{equation}
\begin{equation}\label{eq6}
\langle \langle z, x\rangle , \langle y, x\rangle\rangle =
\langle\langle z, \langle x, y\rangle \rangle , x\rangle +
\langle\langle\langle z, y\rangle , x\rangle , x\rangle -
\langle\langle\langle z, x\rangle , x\rangle , y\rangle,
\end{equation}
where $x$, $y$, $z\in J$.
\end{description}
\end{definition}

A generalized Malcev algebra $M$ is also denoted by $( M,\, +,\, \langle \,, \,\rangle )$, where 
$\langle \,, \,\rangle $ is called the {\bf angle bracket}. Clearly, if the angle bracket 
$\langle \,, \,\rangle $ is anti-commutative, then each of (\ref{eq4}), (\ref{eq5}) and (\ref{eq6}) is equivalent to the following identity:
\begin{equation}\label{eq7}
\langle \langle x, y\rangle , \langle x, z\rangle\rangle =
\langle\langle \langle x, y\rangle , z\rangle , x\rangle +
\langle\langle\langle y, z\rangle , x\rangle , x\rangle +
\langle\langle\langle z, x\rangle , x\rangle , y\rangle .
\end{equation}
Since Malcev algebras are defined by using the anti-commutative law and the identity (\ref{eq7}), the generalized Malcev algebras are indeed a generalization of Malcev algebras. In order to prove that the generalized Malcev algebras are also a generalization of (right) Leibniz algebras, we need to use the {\bf right Jacobian} $J\langle x, y, z\rangle$ which is defined by
\begin{equation}\label{eq8}
J\langle x, y, z\rangle =
\langle\langle x, y\rangle , z\rangle -
\langle x, \langle y, z\rangle \rangle - \langle\langle x, z\rangle , y\rangle  .
\end{equation}

\begin{proposition}\label{pr2.1} Let $(M,\, +,\, \langle \,, \,\rangle  )$ be a right anti-commutative algebra.
\begin{description}
\item[(i)] (\ref{eq4}) is equivalent to
\begin{equation}\label{eq9}
J\langle \langle x, z\rangle , x, y \rangle +J\langle x, y , \langle x, z\rangle \rangle =
\langle J\langle x, y, z\rangle , x \rangle - \langle x, J\langle x, y, z\rangle \rangle .
\end{equation}
\item[(ii)] (\ref{eq5}) is equivalent to
\begin{equation}\label{eq10}
J\langle y, \langle x, z\rangle , x\rangle  = -\langle J\langle y, x, z\rangle , x \rangle .
\end{equation}
\item[(iii)] (\ref{eq6}) is equivalent to
\begin{equation}\label{eq11}
J\langle \langle z, x\rangle , y, x\rangle  = \langle J\langle z, x, y\rangle , x \rangle .
\end{equation}
\end{description}
\end{proposition}

\medskip
\noindent
{\bf Proof} It is clear. 

\hfill\raisebox{1mm}{\framebox[2mm]{}}

\bigskip
An immediate consequence of Proposition 2.1 is that a (right) Leibniz algebra is a generalized Malcev algebra.

\medskip
The next proposition establishes the passage from an alternative $\mathcal{Z}_2$-algebra to a generalized Malcev algebra, which is a generalization of Proposition 1.1 in \cite{hl1}. 

\begin{proposition}\label{pr2.2} If $A=A_0\oplus A_1$ is an alternative $\mathcal{Z}_2$-algebra, then
$( A,\, +,\, \langle \,, \,\rangle)$ is a generalized Malcev algebra, where the angle bracket 
$\langle \, , \, \rangle$ is defined by
\begin{equation}\label{eq12}
\langle x, y\rangle: =xy_0-y_0x
\end{equation}
for $x\in A$, $y=y_0+y_1\in A$ and $y_i\in A_i$ with $i=0$ and $1$.
\end{proposition}

\medskip
\noindent
{\bf Proof} This proposition follows from (\ref{eq12}) and the properties of alternative algebras.

\hfill\raisebox{1mm}{\framebox[2mm]{}}

\bigskip
Let $(M,\, +,\, \langle \,, \,\rangle  )$ be a generalized Malcev algebra over a field $\mathbf{k}$. A subspace $I$ of $M$ is called an {\bf ideal} of $M$ if 
$\langle I, M\rangle\subseteq I$ and $\langle M, I\rangle\subseteq I$. Except $\{0\}$ and $M$, a generalized Malcev algebra $M$ always has another ideal $M^{ann}$, which is defined by
\begin{equation}\label{eq13}
M^{ann}:=\sum _{x, y\in M}\mathbf{k}J\langle x, x, y\rangle + 
\sum _{x\in M}\mathbf{k}\langle x, x\rangle.
\end{equation}
$M^{ann}$ is called the {\bf annihilator} of the generalized Malcev algebra $M$.

\medskip
The building blocks in the context of generalized Malcev algebras are defined by using annihilators.

\bigskip
\begin{definition}\label{def2.2} A generalized Malcev algebra $M$ is said to be {\bf simple} if $M^{ann}\ne 0$ and $M$ has no ideals which are not equal to $\{0\}$, $M^{ann}$ and $M$.
\end{definition}

\medskip
The next proposition gives a basic property of annihilators.

\begin{proposition}\label{pr2.3} If $(M,\, +,\, \langle \,, \,\rangle  )$ is a generalized Malcev algebra over a field $\mathbf{k}$, then the annihilator $M^{ann}$ becomes a module over the Malcev algebra $\displaystyle\frac{M}{M^{ann}}$ under the module action:
\begin{equation}\label{eq14}
\bar{x} a=-a \bar{x}:=-\langle a, x\rangle ,
\end{equation}
where $\bar{x}=x+ M^{ann}\in \displaystyle\frac{M}{M^{ann}}$, $x\in M$ and $a\in M^{ann}$.
\end{proposition}

\medskip
\noindent
{\bf Proof} This is a consequence of the Hu-Liu angle identities. See page 175 of \cite{HCM} for the definition of modules over Malcev algebras.

\hfill\raisebox{1mm}{\framebox[2mm]{}}

\bigskip
\noindent
{\bf 3. Conclusion}

\bigskip
\noindent
Associative $\mathcal{Z}_2$-algebras are a class of associative algebras which can be used to establish quickly the counterparts of the passages from associative algebras to Lie algebras and Jordan algebras. However, using only associative $\mathcal{Z}_2$-algebras is not enough in order to generalize the theory of Lie algebras and Jordan algebras. In fact, another class of associative algebras seems to be more suitable in the study of generalizing universal enveloping algebras and exceptional Jordan algebras. The basic facts about the class of associative algebras will be given in my future work.

\end{document}